\documentclass[12pt]{amsart}
\usepackage{amsmath}
\usepackage{amssymb}
\usepackage{amsthm}
\usepackage{booktabs}
\usepackage[usenames,dvipsnames]{color}
\usepackage{fullpage}
\usepackage{mathtools}
\usepackage{thmtools}
\usepackage{tikz}
\usepackage{url}

\usetikzlibrary{arrows,shapes,automata,backgrounds,decorations,petri,positioning,topaths}

\usepackage[pagebackref]{hyperref} %hyperref wants to be loaded last
\usepackage[all]{hypcap} %fix hyperref links to tables/figures/etc.
\definecolor{darkgreen}{rgb}{0,0.7,0}
\definecolor{purplish}{rgb}{0.5,0,0.8}
\hypersetup{
  breaklinks,colorlinks,citecolor=darkgreen,linkcolor=purplish, % more discreet colours for links
  pdftitle={The generating function for total displacement},
  pdfauthor={T.~Kyle Petersen and Mathieu Guay-Paquet},
}
\theoremstyle{plain}
\newtheorem{thm}{Theorem}[section]

\newtheorem{prop}[thm]{Proposition}

\theoremstyle{definition}

\newtheorem{remark}[thm]{Remark}
\newtheorem{ex}[thm]{Example}

\numberwithin{equation}{section} % requires package amsthm

\DeclarePairedDelimiter{\abs}{|}{|}
\DeclarePairedDelimiter{\set}{\{}{\}}

\DeclareMathOperator{\Motz}{Motz}

\DeclareMathOperator{\dep}{dep}

\DeclareMathOperator{\ar}{area}
\DeclareMathOperator{\wt}{\omega}

\title{The generating function for total displacement}

\author[T.~K. Petersen]{T.~Kyle Petersen}
\email{tpeter21@depaul.edu}
\address{
Department of Mathematical Sciences \\
DePaul University \\
Chicago IL 60614 \\
USA}

\author{Mathieu Guay-Paquet}
\email{mathieu.guaypaquet@lacim.ca}
\address{
LaCIM \\
Universit\'e du Qu\'ebec \`a Montr\'eal \\
201 Pr\'esident-Kennedy \\
Montr\'eal QC\ \ H2X~3Y7 \\
Canada}


\begin{document}

\begin{abstract}
In a 1977 paper, Diaconis and Graham studied what Knuth calls the \emph{total displacement} of a permutation $w$, which is the sum of the distances $\abs{w(i)-i}$. In recent work of the first author and Tenner, this statistic appears as twice the type $A_{n-1}$ version of a statistic for Coxeter groups called the \emph{depth} of $w$. There are various enumerative results for this statistic in the work of Diaconis and Graham, codified as exercises in Knuth's textbook, and some other results in the work of Petersen and Tenner. However, no formula for the generating function of this statistic appears in the literature. Knuth comments that ``the generating function for total displacement does not appear to have a simple form." In this paper, we translate the problem of computing the distribution of total displacement into a problem of counting weighted Motzkin paths. In this way, standard techniques allow us to express the generating function for total displacement as a continued fraction.
\end{abstract}

\maketitle

%---------------------------------------------------------------
\section{Introduction}
%---------------------------------------------------------------

Let $S_n$ denote the symmetric group of permutations of $\set{1, 2, \ldots, n}$. Diaconis and Graham~\cite{DG} studied what they called \emph{Spearman's disarray} for permutations in $S_n$. The disarray statistic, which Knuth calls \emph{total displacement}~\cite[Problem 5.1.1.28]{K}, is defined for a permutation $w$ as \[\sum_{i=1}^n \abs{w(i)-i} = 2\sum_{w(i) > i} (w(i)-i).\]

In work of the first author and Tenner~\cite{PT}, half of the total displacement is shown to be equal to the \emph{depth} of a permutation, which is a measure of certain minimal factorizations of a permutation into transpositions. That is, define
\begin{equation}\label{eq:depdef}
  \dep(w) = \min\set*{\sum_{s=1}^k (j_s-i_s) : w = (i_1\, j_1) \cdots (i_k\, j_k)}.
\end{equation}
Then~\cite[Theorem 1.1]{PT} shows that
\begin{equation}\label{eq:dep}
  \dep(w) = \sum_{\mathclap{w(i) > i}} (w(i)-i).
\end{equation}
The definition of depth is easily extended to all finite Coxeter groups (as discussed in~\cite{PT}, one replaces transpositions with reflections); this is the type $A_{n-1}$ version.

For example, if $w = 3715246$ in one-line notation, we can write $w$ as a product of transpositions as follows: \[w = (67)(46)(24)(13)(45).\] Moreover,~\cite[Section 3.1]{PT} shows that this minimizes the right-hand side of \autoref{eq:depdef}, so \[\dep(w) = (7-6) + (6-4) + (4-2) + (3-1) + (5-4) = 8.\] But according to \autoref{eq:dep}, we can compute the same number from the one-line notation: \[\dep(w) = (3-1) + (7-2) + (5-4) = 8.\]

The goal of this paper is to construct the generating function for the distribution of depth, or equivalently, total displacement. In Knuth's remarks on the problem, ``The generating function for total displacement does not appear to have a simple form"~\cite[Solution for 5.1.1.28]{K}. The generating function we produce is a continued fraction, so whether one disagrees with Knuth depends on whether one considers a continued fraction to be ``simple."

We include here for reference \autoref{tab:hnk}, which contains the distribution of depths for $n \leq 8$. This array can be found as entry \href{http://oeis.org/A062869}{A062869} of~\cite{oeis}.

\begin{table}[h]
\begin{center}
{\tiny\begin{tabular}{r|*{17}{c}}
\toprule
& $k=0$ & 1 & 2 & 3 &  4 & 5 & 6 & 7 & 8 & 9 & 10 & 11 & 12 & 13 & 14 & 15 & 16 \\
\midrule
$n=1$ & 1 \\
   2  & 1 & 1 \\
   3  & 1 & 2 &  3 \\
   4  & 1 & 3 &  7 &   9 &   4 \\
   5  & 1 & 4 & 12 &  24 &  35 &  24 &   20 \\
   6  & 1 & 5 & 18 &  46 &  93 & 137 &  148 &  136 &  100 &   36 \\
   7  & 1 & 6 & 25 &  76 & 187 & 366 &  591 &  744 &  884 &  832 &  716 &  360 &  252 \\
   8  & 1 & 7 & 33 & 115 & 327 & 765 & 1523 & 2553 & 3696 & 4852 & 5708 & 5892 & 5452 & 4212 & 2844 & 1764 & 576 \\
\bottomrule
\end{tabular}}
\end{center}
\caption{The number of permutations $w \in S_n$ with depth $\dep(w) = k$.}\label{tab:hnk}
\end{table}

%---------------------------------------------------------------
\section{A map from permutations to Motzkin paths}\label{sec:map-def}
%---------------------------------------------------------------

We begin by translating the problem of counting permutations by depth to the problem of counting weighted Motzkin paths by area. More specifically, we will give a map from permutations to Motzkin paths which takes depth to area, and then weigh each Motzkin path by the number of permutations which map to it.

Recall that a \emph{Motzkin path} of length $n$ is a sequence of $n$ letters from the set $\set{U, D, H}$ such that the subword on the letters $U$ and $D$ forms a balanced parenthesization. We can draw Motzkin paths by identifying the letters $U$, $D$, and $H$ with the segments
\begin{tikzpicture}[scale=1.2,baseline=0] \draw (0,0) node {$\bullet$}; \draw (1ex,1ex) node {$\bullet$}; \draw (0,0) -- (1ex,1ex); \end{tikzpicture},
\begin{tikzpicture}[scale=1.2,baseline=0] \draw (0,1ex) node {$\bullet$}; \draw (1ex,0) node {$\bullet$}; \draw (0,1ex) -- (1ex,0); \end{tikzpicture}, and
\begin{tikzpicture}[scale=1.8,baseline=-.5ex] \draw (0,0) node {$\bullet$}; \draw (1ex,0) node {$\bullet$}; \draw (0,0) -- (1ex,0); \end{tikzpicture},
respectively. Then, the balanced parenthesization condition means only that the path starts and ends at the same height, and never goes below this height. The \emph{area} of a Motzkin path is the area of the region above the base line and below the path. For example, if $p = UUHUDDHDH$, we would draw the following picture and see that $\ar(p) = 12$:
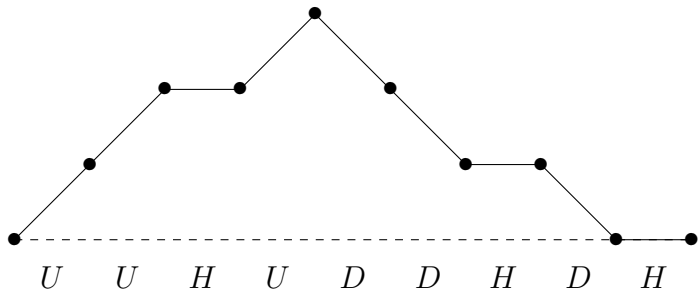
\[
\begin{tikzpicture}
\draw (-.5,-.5)
++(1,0) node {$U$}
++(1,0) node {$U$}
++(1,0) node {$H$}
++(1,0) node {$U$}
++(1,0) node {$D$}
++(1,0) node {$D$}
++(1,0) node {$H$}
++(1,0) node {$D$}
++(1,0) node {$H$};
\draw (0,0) node {$\bullet$} --
++(1,1) node {$\bullet$} --
++(1,1) node {$\bullet$} --
++(1,0) node {$\bullet$} --
++(1,1) node {$\bullet$} --
++(1,-1) node {$\bullet$} --
++(1,-1) node {$\bullet$} --
++(1,0) node {$\bullet$} --
++(1,-1) node {$\bullet$} --
++(1,0) node {$\bullet$};
\draw[dashed] (0,0) -- (9,0);
\end{tikzpicture}
\]

The map we are looking for is $\phi \colon S_n \to \Motz_n$ by $\phi(w) = p_1 \cdots p_n$, where
\[
  p_i = \begin{cases}
    U & \text{if $w^{-1}(i) > i < w(i)$,} \\
    D & \text{if $w^{-1}(i) < i > w(i)$,} \\
    H & \text{otherwise.}
  \end{cases}
\]
For example, if $w = 3715246$, we have $\phi(w) = UUDUDHD$. We can see the correspondence more easily if we draw our permutations as follows. Write the numbers $1, \ldots, n$ on a line, and draw an arrow from $i$ to $j$ if $w(i) = j$ and $i \neq j$. Moreover, if $i < j$, draw the arrow above the numbers, and if $i > j$, draw the arrow below. Continuing with the same example of $w=3715246$, we draw the following:
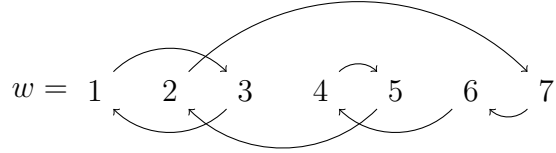
\[
w = \begin{tikzpicture}[baseline=-.5ex,bend angle=45]
\foreach \x in {1,...,7}{\draw (\x,0) node (\x) {\x};}
\draw[->] (1) to[bend left] (3);
\draw[->] (2) to[bend left] (7);
\draw[->] (4) to[bend left] (5);
\draw[->] (7) to[bend left] (6);
\draw[->] (6) to[bend left] (4);
\draw[->] (5) to[bend left] (2);
\draw[->] (3) to[bend left] (1);
\end{tikzpicture}
\]
In terms of the picture, we have $p_i = U$ if the element $i$ has an incoming arrow from below, and an outgoing arrow above (both to its right); we have $p_i = D$ if the element $i$ has an incoming arrow from above, and an outgoing arrow below (both to its left); and we have $p_i = H$ if $i$ is a fixed point or if the arrows incident to it are both above or both below.

\autoref{tab:S4} shows the correspondence for all permutations in $S_4$, grouped by their associated Motzkin paths.

%---------------------------------------------------------------
\begin{table}[b]
\begin{center}
\begin{tabular}{ccc}
\toprule
$w \in S_4$ & $\phi(w) \in \Motz_4$  & $\dep(w) = \ar(\phi(w))$ \\
\midrule
$\begin{array}{c} 1234 \end{array}$ &
\begin{tikzpicture}[baseline=-.5ex]
  \draw (0,0) node {$\bullet$} --
  ++(1,0) node {$\bullet$} --
  ++(1,0) node {$\bullet$} --
  ++(1,0) node {$\bullet$} --
  ++(1,0) node {$\bullet$};
  \draw[dashed] (0,0) -- (4,0);
\end{tikzpicture} &
0 \\
\midrule
$\begin{array}{c} 1243 \end{array}$ &
\begin{tikzpicture}[baseline=.5cm-.5ex]
  \draw (0,0) node {$\bullet$} --
  ++(1,0) node {$\bullet$} --
  ++(1,0) node {$\bullet$} --
  ++(1,1) node {$\bullet$} --
  ++(1,-1) node {$\bullet$};
  \draw[dashed] (0,0) -- (4,0);
\end{tikzpicture} &
1 \\
\midrule
$\begin{array}{c} 1324 \end{array}$ &
\begin{tikzpicture}[baseline=.5cm-.5ex]
  \draw (0,0) node {$\bullet$} --
  ++(1,0) node {$\bullet$} --
  ++(1,1) node {$\bullet$} --
  ++(1,-1) node {$\bullet$} --
  ++(1,0) node {$\bullet$};
  \draw[dashed] (0,0) -- (4,0);
\end{tikzpicture} &
1 \\
\midrule
$\begin{array}{c} 2134 \end{array}$ &
\begin{tikzpicture}[baseline=.5cm-.5ex]
  \draw (0,0) node {$\bullet$} --
  ++(1,1) node {$\bullet$} --
  ++(1,-1) node {$\bullet$} --
  ++(1,0) node {$\bullet$} --
  ++(1,0) node {$\bullet$};
  \draw[dashed] (0,0) -- (4,0);
\end{tikzpicture} &
1 \\
\midrule
$\begin{array}{c} 2143 \end{array}$ &
\begin{tikzpicture}[baseline=.5cm-.5ex]
  \draw (0,0) node {$\bullet$} --
  ++(1,1) node {$\bullet$} --
  ++(1,-1) node {$\bullet$} --
  ++(1,1) node {$\bullet$} --
  ++(1,-1) node {$\bullet$};
  \draw[dashed] (0,0) -- (4,0);
\end{tikzpicture} &
2 \\
\midrule
$\begin{array}{c} 1342 \\ 1423 \\ 1432 \end{array}$ &
\begin{tikzpicture}[baseline=.5cm-.5ex]
  \draw (0,0) node {$\bullet$} --
  ++(1,0) node {$\bullet$} --
  ++(1,1) node {$\bullet$} --
  ++(1,0) node {$\bullet$} --
  ++(1,-1) node {$\bullet$};
  \draw[dashed] (0,0) -- (4,0);
\end{tikzpicture} &
2 \\
\midrule
$\begin{array}{c} 2314 \\ 3124 \\ 3214 \end{array}$ &
\begin{tikzpicture}[baseline=.5cm-.5ex]
  \draw (0,0) node {$\bullet$} --
  ++(1,1) node {$\bullet$} --
  ++(1,0) node {$\bullet$} --
  ++(1,-1) node {$\bullet$} --
  ++(1,0) node {$\bullet$};
  \draw[dashed] (0,0) -- (4,0);
\end{tikzpicture} &
2 \\
\midrule
$\begin{array}{c} 2341 \\ 2413 \\ 2431 \\ 3142 \\ 3241 \\ 4123 \\ 4132 \\ 4213 \\ 4231 \end{array}$ &
\begin{tikzpicture}[baseline=.5cm-.5ex]
  \draw (0,0) node {$\bullet$} --
  ++(1,1) node {$\bullet$} --
  ++(1,0) node {$\bullet$} --
  ++(1,0) node {$\bullet$} --
  ++(1,-1) node {$\bullet$};
  \draw[dashed] (0,0) -- (4,0);
\end{tikzpicture} &
3 \\
\midrule
$\begin{array}{c} 3412 \\ 3421 \\ 4312 \\ 4321 \end{array}$ &
\begin{tikzpicture}[baseline=1cm-.5ex]
  \draw (0,0) node {$\bullet$} --
  ++(1,1) node {$\bullet$} --
  ++(1,1) node {$\bullet$} --
  ++(1,-1) node {$\bullet$} --
  ++(1,-1) node {$\bullet$};
  \draw[dashed] (0,0) -- (4,0);
\end{tikzpicture} &
4 \\
\bottomrule
\end{tabular}
\end{center}
\caption{The map $\phi$ on $S_4$.}\label{tab:S4}
\end{table}

%---------------------------------------------------------------
\section{Area equals depth}\label{sec:area-depth}
%---------------------------------------------------------------

Let us show that the map $\phi \colon S_n \to \Motz_n$ takes depth to area.

An \emph{excedance} of a permutation $w$ is a position $i$ such that $w(i) > i$. The characterization of depth in \autoref{eq:dep} shows that depth is equal to the sum of the sizes $w(i)-i$ of the excedances. If we draw the permutation $w$ as in \autoref{sec:map-def}, an excedance is simply an arrow from a smaller number to a larger one, that is, an arrow drawn above the line of numbers. If $i$ is an excedance, so that $w(i) = j > i$, this is drawn as:
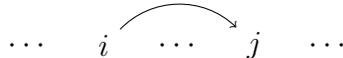
\[
\begin{tikzpicture}[baseline=0,bend angle=45]
\draw (0,0) node {$\cdots$}
++(1,0) node (i) {$i$}
++(1,0) node {$\cdots$}
++(1,0) node (j) {$j$}
++(1,0) node {$\cdots$};
\draw[->] (i) to[bend left] (j);
\end{tikzpicture}
\]
and the contribution of the excedance $i$ to the depth of $w$ is the length of the arrow from $i$ to $j$. If $j$ is itself an excedance, so that $w(j) = k > j$, the drawing looks like:
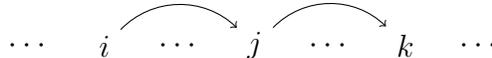
\[
\begin{tikzpicture}[baseline=0,bend angle=45]
\draw (0,0) node {$\cdots$}
++(1,0) node (i) {$i$}
++(1,0) node {$\cdots$}
++(1,0) node (j) {$j$}
++(1,0) node {$\cdots$}
++(1,0) node (k) {$k$}
++(1,0) node {$\cdots$};
\draw[->] (i) to[bend left] (j);
\draw[->] (j) to[bend left] (k);
\end{tikzpicture}
\]
and the total contribution to depth is $(j-i) + (k-j)= k-i$; that is, the distance from $i$ to $k$ is all that matters. In fact, we can compute the depth of $w$ if we only know where each string of right arrows begins and ends. And if $\phi(w) = p_1 \cdots p_n$ is the Motzkin path associated to $w$, then $i$ is the beginning of a string of right arrows exactly when $w^{-1}(i) > i < w(i)$, i.e., $p_i = U$, and $k$ is the end of a string of right arrows exactly when $w^{-1}(k) < k > w(k)$, i.e., $p_k = D$. Thus, the depth can be characterized as follows.

\begin{prop}\label{prop:depthDU}
Let $w \in S_n$ be a permutation and $\phi(w) = p_1 \cdots p_n$ be the associated Motzkin path. Then, \[\dep(w) = \sum_{p_k = D} k - \sum_{p_i = U} i.\]
\end{prop}

For example, if $w = 3542176$, we draw:
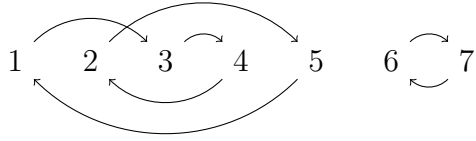
\[
\begin{tikzpicture}[baseline=0,bend angle=45]
\foreach \x in {1,...,7}{
\draw (\x,0) node (\x) {\x};
}
\draw[->] (1) to[bend left] (3);
\draw[->] (3) to[bend left] (4);
\draw[->] (2) to[bend left] (5);
\draw[->] (6) to[bend left] (7);
\draw[->] (4) to[bend left] (2);
\draw[->] (5) to[bend left] (1);
\draw[->] (7) to[bend left] (6);
\end{tikzpicture}
\]
and the corresponding Motzkin path is $p = \phi(w) = UUHDDUD$, so the depth is \[\dep(w) = (4+5+7) - (1+2+6).\] It turns out that this number is equal to the area under the path $p$, as can be seen by rearranging the terms as \[\dep(w) = (5-1) + (4-2) + (7-6)\]
and consulting the following diagram:
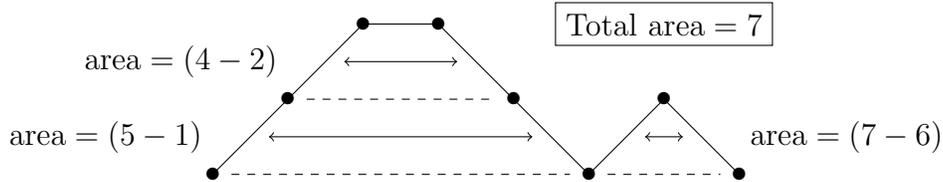
\[
\begin{tikzpicture}[baseline=0]
\draw (0,0) node (L1) {$\bullet$} --
++(1,1) node (L2) {$\bullet$} --
++(1,1) node {$\bullet$} --
++(1,0) node {$\bullet$} --
++(1,-1) node (R2) {$\bullet$} --
++(1,-1) node (M1) {$\bullet$} --
++(1,1) node {$\bullet$} --
++(1,-1) node (R1) {$\bullet$};
\draw[dashed] (L1) -- (M1) -- (R1);
\draw[dashed] (L2) -- (R2);
\draw[<->] (.75,.5)--(4.25,.5);
\draw[<->] (1.75,1.5)--(3.25,1.5);
\draw[<->] (5.75,.5)--(6.25,.5);
\draw (0,.5) node[left] {$\ar = (5-1)$};
\draw (1,1.5) node[left] {$\ar = (4-2)$};
\draw (7,.5) node[right] {$\ar = (7-6)$};
\draw (6,2) node[rectangle,draw=black] {Total $\ar = 7$};
\end{tikzpicture}
\]

In general, for a Motzkin path $p = p_1 \cdots p_n$, the letters $U$ and $D$ in $p$ form a balanced parenthesization, so they can be naturally paired up; furthermore, the pair of steps $p_i = U$ and $p_k = D$ contribute a horizontal strip of area $k-i$ below the Motzkin path. Thus, we can characterize the area under $p$ as \[\ar(p) = \sum_{p_k = D} k - \sum_{p_i = U} i.\] Together with \autoref{prop:depthDU}, this observation proves that the map $\phi$ takes depth to area.

\begin{prop}\label{prop:area-depth}
For any $w \in S_n$, $\dep(w) = \ar(\phi(w))$.
\end{prop}

\begin{remark}
  In~\cite[Proposition~3.2]{PT}, the first author and Tenner prove that for all $w \in S_n$, the depth of $w$ is at most $\lfloor n^2/4\rfloor$, and this bound is sharp. We can recover this result as a corollary of \autoref{prop:area-depth}, since $\lfloor n^2/4\rfloor$ is the largest possible area for a Motzkin path with $n$ steps, attained by the path $p = U^k D^k$ if $n = 2k$ is even, and by the path $p = U^k H D^k$ if $n = 2k+1$ is odd.
\end{remark}

%---------------------------------------------------------------
\section{How many permutations map to a Motzkin path?}\label{sec:preimage}
%---------------------------------------------------------------

Having shown that the map $\phi \colon S_n \to \Motz_n$ encodes the depth of a permutation as the area of the associated Motzkin path, the next step in determining the distribution of the depth statistic is to compute the size of the preimage $\phi^{-1}(p)$ for each Motzkin path $p$.

Fortunately, this can be done easily, as demonstrated with the following example.

\begin{ex}\label{ex:preimage}
Given the Motzkin path $p = UUHDDUD$, we start by writing the numbers $1, \ldots, 7$ on a line and drawing incoming and outgoing half-arcs for every $U$ and $D$ (without connecting the half-arcs to each other):
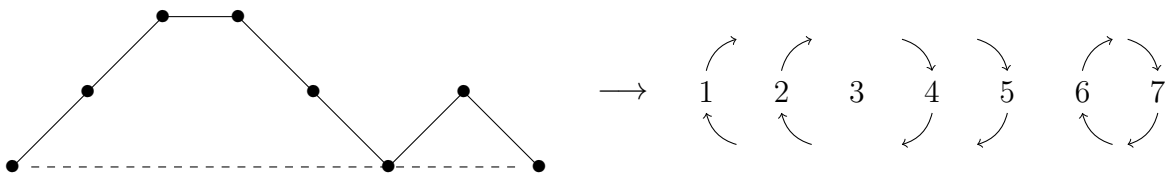
\[
\begin{tikzpicture}[baseline=1cm-.5ex]
\draw (0,0) node (left) {$\bullet$} --
++(1,1) node {$\bullet$} --
++(1,1) node {$\bullet$} --
++(1,0) node {$\bullet$} --
++(1,-1) node {$\bullet$} --
++(1,-1) node {$\bullet$} --
++(1,1) node {$\bullet$} --
++(1,-1) node (right) {$\bullet$};
\draw[dashed] (left) -- (right);
\end{tikzpicture}
\quad\longrightarrow\quad
\begin{tikzpicture}[baseline=-.5ex,bend angle=30]
\foreach \x in {1,...,7}{
\draw (\x,0) node (\x) {\x};
}
\draw[->] (1) to[bend left]  +( .4, .7);
\draw[<-] (1) to[bend right] +( .4,-.7);
\draw[->] (2) to[bend left]  +( .4, .7);
\draw[<-] (2) to[bend right] +( .4,-.7);
\draw[->] (6) to[bend left]  +( .4, .7);
\draw[<-] (6) to[bend right] +( .4,-.7);
\draw[<-] (4) to[bend right] +(-.4, .7);
\draw[->] (4) to[bend left]  +(-.4,-.7);
\draw[<-] (5) to[bend right] +(-.4, .7);
\draw[->] (5) to[bend left]  +(-.4,-.7);
\draw[<-] (7) to[bend right] +(-.4, .7);
\draw[->] (7) to[bend left]  +(-.4,-.7);
\end{tikzpicture}
\]
As noted in \autoref{sec:area-depth}, when drawing a permutation $w$ such that $\phi(w) = p$, each excedance of $w$ is drawn as a right-pointing arrow above the line of numbers, and these arrows can be grouped into strings of right-pointing arrows, which start at a position $i$ with $p_i = U$ and end at a position $k$ with $p_k = D$, possibly with intermediate steps at positions $j$ with $p_j = H$. So, to form the permutations that correspond to this path, we will first ignore the positions $j$ with $p_j = H$ and match up all the half-arcs above the line of numbers, outgoing with incoming, to indicate the starting and ending positions of each string of right-pointing arrows. In this example, we are forced to form the pair $6 \to 7$, but we have two choices for pairing the outgoing half-arcs $1 \to \cdot$, $2 \to \cdot$ and the incoming half-arcs $\cdot \to 4$, $\cdot \to 5$; that is:
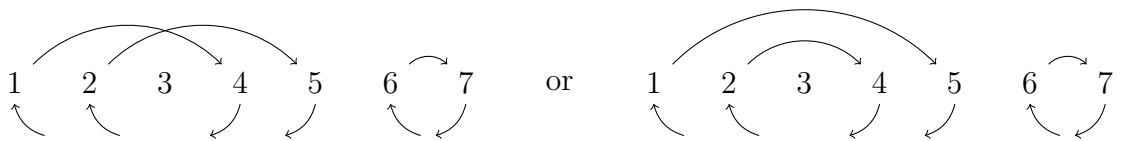
\[
\begin{tikzpicture}[baseline=-.5ex,bend angle=30]
\foreach \x in {1,...,7}{
\draw (\x,0) node (\x) {\x};
}
\draw[->] (1) to[bend left=45] (4);
\draw[->] (2) to[bend left=45] (5);
\draw[->] (6) to[bend left=45] (7);
\draw[<-] (1) to[bend right] +( .4,-.7);
\draw[<-] (2) to[bend right] +( .4,-.7);
\draw[<-] (6) to[bend right] +( .4,-.7);
\draw[->] (4) to[bend left]  +(-.4,-.7);
\draw[->] (5) to[bend left]  +(-.4,-.7);
\draw[->] (7) to[bend left]  +(-.4,-.7);
\end{tikzpicture}
\qquad\text{or}\qquad
\begin{tikzpicture}[baseline=-.5ex,bend angle=30]
\foreach \x in {1,...,7}{
\draw (\x,0) node (\x) {\x};
}
\draw[->] (1) to[bend left=45] (5);
\draw[->] (2) to[bend left=45] (4);
\draw[->] (6) to[bend left=45] (7);
\draw[<-] (1) to[bend right] +( .4,-.7);
\draw[<-] (2) to[bend right] +( .4,-.7);
\draw[<-] (6) to[bend right] +( .4,-.7);
\draw[->] (4) to[bend left]  +(-.4,-.7);
\draw[->] (5) to[bend left]  +(-.4,-.7);
\draw[->] (7) to[bend left]  +(-.4,-.7);
\end{tikzpicture}.
\]
Independently, we also have two choices for matching up the half-arcs below the line of numbers into (for now, single-step) strings of left-pointing arrows in the permutation diagram for $w$. Finally, to complete the diagram, we only need to decide what to do with the position 3, for which $p_3 = H$: it can either be a fixed point of $w$, or join one of the two strings of right-pointing arrows above it, or join one of the two strings of left-pointing arrows below it, for a total of five choices. (Note that for all matchings here, there are two strings of right-pointing arrows above position 3 and two strings of left-pointing arrows below it.) All in all, we have $2 \cdot 2 \cdot 5 = 20$ possible diagrams (all valid) for a permutation $w$ such that $\phi(w) = p = UUHDDUD$.
\end{ex}

In general, we can count the number of permutations corresponding to a given Motzkin path by reconstructing the possible permutation diagrams as in \autoref{ex:preimage}; first we count the number of ways of matching up the outgoing and incoming half-arcs above the line of numbers, then we count the number of ways of matching up the half-arcs below, and finally we count the number of ways of dealing with the positions corresponding to $H$ steps.

To do this, we will need to define the \emph{height} of each step in a Motzkin path. We draw our paths starting at $(0,0)$ and define the height $h_i$ of a step $p_i$ to be the maximum height achieved on that part of the path. That is, $h_i = j$ if
\begin{itemize}
\item $p_i = U$ from $(i-1,j-1)$ to $(i,j)$,
\item $p_i = H$ from $(i-1,j)$ to $(i,j)$, or
\item $p_i = D$ from $(i-1,j)$ to $(i,j-1)$.
\end{itemize}
For example, the steps of $p=UUHDDUD$ have heights $(h_1,\ldots,h_7) = (1,2,2,2,1,1,1)$.

Note that the height $h_i$ of a step in $p$ is also the number of strings of right-pointing arrows which appear above position $i$ in the diagram of any permutation $w$ with $\phi(w) = p$. Indeed, if we scan the diagram from left to right, the height increases by one for each $U$ step, which marks the beginning of a string of right-pointing arrows, stays constant for each $H$ step, and decreases by one for each $D$ step, which marks the end of a string of right-pointing arrows. This observation lets us count the number of ways of matching up the half-arcs above the line of numbers; if we build the permutation diagram from left to right, then we have no choices to make for $U$ and $H$ steps, but for a step $p_i = D$, we can freely choose any of the $h_i$ strings of right-pointing arrows above position $i$ to terminate there. In other words, the number of ways to match up the half-arcs above the line of numbers is \[\prod_{p_i = D} h_i.\]

We can apply the same argument to the half-arcs below the line of numbers. If we construct a matching from right to left, the number of possibilities is simply \[\prod_{p_i = U} h_i.\]

Finally, for each step $p_i = H$, there are $2h_i + 1$ choices for dealing with position $i$ in the diagram: let $i$ be a fixed point of the permutation, or join one of the $h_i$ strings of right-pointing arrows above it, or join one of the $h_i$ strings of left-pointing arrows below it.

Hence, we define the \emph{weight} of step $p_i$ to be
\[
  \wt_i = \begin{cases}
    h_i &\text{if $p_i = U$ or $p_i = D$,} \\
    2h_i + 1 &\text{if $p_i = H$,}
  \end{cases}
\]
and the weight of a path $p = p_1 \cdots p_n$ to be the product
\[
\wt(p) = \wt_1 \cdots \wt_n.
\]
Then, the weight of a Motzkin path $p$ is the number of permutations in its preimage $\phi^{-1}(p)$.

\begin{prop}\label{prop:preimage}
Let $p \in \Motz_n$. Then \[\wt(p) = \abs{\set{w \in S_n : \phi(w) = p}}.\]
\end{prop}

In the example $p=UUHDDUD$, we have $\wt(p) = 1 \cdot 2 \cdot 5 \cdot 2 \cdot 1 \cdot 1 \cdot 1=20$, which can be seen visually as:
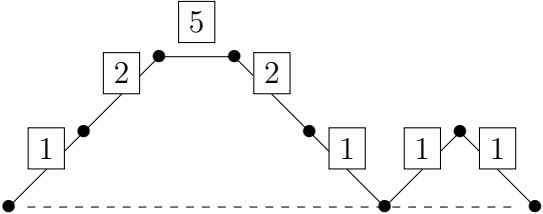
\[
\begin{tikzpicture}
[baseline=1cm-.5ex,
 weight/.style={above,rectangle,draw,fill=white}]
\draw (0,0) node (start) {$\bullet$} --
node[weight] {1} ++(1,1) node {$\bullet$} --
node[weight] {2} ++(1,1) node {$\bullet$} --
node[weight,above=1ex] {5} ++(1,0) node {$\bullet$} --
node[weight] {2} ++(1,-1) node {$\bullet$} --
node[weight] {1} ++(1,-1) node {$\bullet$} --
node[weight] {1} ++(1,1) node {$\bullet$} --
node[weight] {1} ++(1,-1) node (end) {$\bullet$};
\draw[dashed] (start) -- (end);
\end{tikzpicture}
\]
As a larger example, the path $q=UHDHUUUDHUDDHD$ has the step weights:
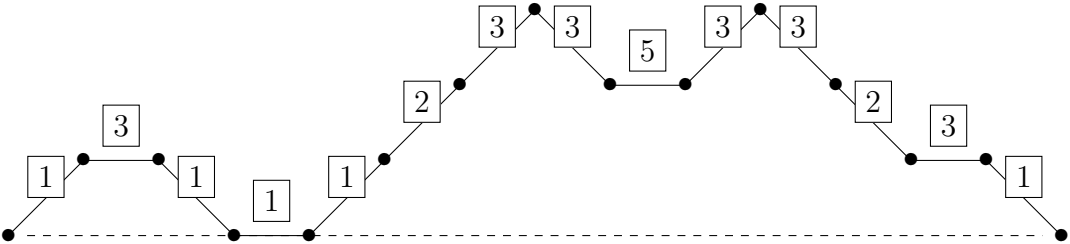
\[
\begin{tikzpicture}
[baseline=1cm-.5ex,
 weight/.style={above,rectangle,draw,fill=white}]
\draw (0,0) node (start) {$\bullet$} --
node[weight] {1} ++(1,1) node {$\bullet$} --
node[weight,above=1ex] {3} ++(1,0) node {$\bullet$} --
node[weight] {1} ++(1,-1) node {$\bullet$} --
node[weight,above=1ex] {1} ++(1,0) node {$\bullet$} --
node[weight] {1} ++(1,1) node {$\bullet$} --
node[weight] {2} ++(1,1) node {$\bullet$} --
node[weight] {3} ++(1,1) node {$\bullet$} --
node[weight] {3} ++(1,-1) node {$\bullet$} --
node[weight,above=1ex] {5} ++(1,0) node {$\bullet$} --
node[weight] {3} ++(1,1) node {$\bullet$} --
node[weight] {3} ++(1,-1) node {$\bullet$} --
node[weight] {2} ++(1,-1) node {$\bullet$} --
node[weight,above=1ex] {3} ++(1,0) node {$\bullet$} --
node[weight] {1} ++(1,-1) node (end) {$\bullet$};
\draw[dashed] (start) -- (end);
\end{tikzpicture}
\]
so $\wt(q) = 2^2 \cdot 3^6 \cdot 5 = 14580$.

\begin{remark}
  In~\cite[Proposition~3.3]{PT}, the first author and Tenner prove that the number of permutations $w \in S_n$ which achieve the maximal depth of $\lfloor n^2/4\rfloor$ is
  \[
    \abs{\set{w \in S_n : \dep(w) = \lfloor n^2/4\rfloor}} = \begin{cases}
      (k!)^2 &\text{if $n=2k$,} \\
      n(k!)^2 &\text{if $n=2k+1$}.
    \end{cases}
  \]
  We can recover this result as a corollary of \autoref{prop:preimage} by noting that this is the weight of the Motzkin path with maximal area, namely $p = U^k D^k$ if $n = 2k$ is even, and $p = U^k H D^k$ if $n = 2k+1$ is odd.

  Statements equivalent to~\cite[Proposition~3.2]{PT} and~\cite[Proposition~3.3]{PT} can be found in the paper of Diaconis and Graham, although without proof (see Table~1 and Remark~2 of~\cite{DG}). They are also mentioned in the remarks (and links therein) for entry \href{http://oeis.org/A062870}{A062870} of~\cite{oeis}.
\end{remark}

%---------------------------------------------------------------
\section{Counting weighted Motzkin paths by area}
%---------------------------------------------------------------

Taking Propositions~\ref{prop:area-depth} and~\ref{prop:preimage} into account, we can express the generating function for permutations with respect to depth as
\[
  F(t, z)
    = \sum_{n \geq 0} \sum_{w \in S_n} t^{\dep(w)} z^n
    = \sum_{p \in \Motz} \wt(p) t^{\ar(p)} z^{\abs{p}},
\]
where $\abs{p}$ is the number of steps in the path $p$. Furthermore, if we decompose each Motzkin path into vertical strips (instead of horizontal strips as in \autoref{sec:area-depth}) to compute its area, we can rewrite the whole term $\wt(p) t^{\ar(p)} z^{\abs{p}}$ as a product over the steps of $p$. For example, if $p = UUHDDUD$, we would have the modified weights
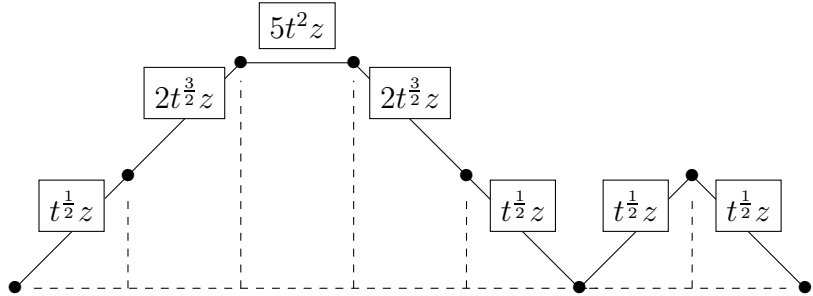
\[
\begin{tikzpicture}
[baseline=1cm-.5ex,
 weight/.style={above,rectangle,draw,fill=white},
 scale=1.5]
\draw (0,0) node (start) {$\bullet$} --
node[weight] {$t^{\frac12} z$} ++(1,1) node (1) {$\bullet$} --
node[weight] {$2t^{\frac32} z$} ++(1,1) node (2) {$\bullet$} --
node[weight,above=1ex] {$5t^2 z$} ++(1,0) node (3) {$\bullet$} --
node[weight] {$2t^{\frac32} z$} ++(1,-1) node (4) {$\bullet$} --
node[weight] {$t^{\frac12} z$} ++(1,-1) node (5) {$\bullet$} --
node[weight] {$t^{\frac12} z$} ++(1,1) node (6) {$\bullet$} --
node[weight] {$t^{\frac12} z$} ++(1,-1) node (end) {$\bullet$};
\draw[dashed] (start) -- (end);
\foreach \x in {1,...,6}{
\draw[dashed] (\x,0) -- (\x);
}
\end{tikzpicture}
\]

Following~\cite[Section~5.2]{GJ}, we can count Motzkin paths with these modified weights using continued fractions. For brevity, let $a_i$, $b_i$, and $c_i$ represent the weight of a $U$, $H$, or $D$ step at height $i$, respectively. Any Motzkin path can be uniquely decomposed as a list of subpaths by cutting it at every (integer) point where it reaches height 0. The pieces in this decomposition will be of the form $H$ or $Up_1D$, where $p_1$ is a Motzkin path with base height 1 instead of 0. Hence, the generating function for weighted Motzkin paths is
\[
  F = \frac{1}{1 - b_0 - a_1 c_1 F_1},
\]
where $F_1$ is the generating function for weighted Motzkin paths with base height 1. We can apply the same reasoning to $F_1$ to express it in terms of the generating function $F_2$ for weighted Motzkin paths with base height 2:
\[
  F = \frac{1}{1 - b_0 - \displaystyle\frac{a_1 c_1}
              {1 - b_1 - a_2 c_2 F_2}}.
\]
Continuing in this manner, we obtain the continued fraction
\[
  F = \frac{1}{1 - b_0 - \displaystyle\frac{a_1 c_1}
              {1 - b_1 - \displaystyle\frac{a_2 c_2}
              {1 - b_2 - \displaystyle\frac{a_3 c_3}
              {1 - \cdots}}}}.
\]
If we substitute our actual modified weights for $a_i$, $b_i$ and $c_i$, we get
\begin{equation}\label{eq:j-fraction}
  F(t, z) = \frac{1}{1 - z - \displaystyle\frac{tz^2}
                    {1 - 3tz - \displaystyle\frac{4t^3 z^2}
                    {1 - 5t^2 z - \displaystyle\frac{9t^5 z^2}
                    {1 - 7t^3 z - \displaystyle\frac{16t^7 z^2}
                    {1 - \cdots}}}}}.
\end{equation}

\begin{remark}
  If we set $t = 1$ in $F(t, z)$, so that we are ignoring depth, we recover a continued fraction due to Euler (see~\cite[Section~5.2.11]{GJ}) for the series \[F(1, z) = \sum_{n\geq 0} n! z^n.\] In~\cite{GJ} this is obtained algebraically, whereas here we get a refined version $F(t, z)$ of this series combinatorially (although see~\cite[Section~5.2.16]{GJ} for a closely related continued fraction with a combinatoral interpretation).
\end{remark}

Note that in the continued fraction~\eqref{eq:j-fraction}, the denominators contain both a linear term and a quadratic term in $z$, so in the language of~\cite{GJ} it would be called a $J$-fraction. However, because of the special form of these terms in our case, we can rewrite it as a simpler $S$-fraction, which contains only linear terms in $z$ in its denominators (see~\cite[Proposition~5.2.2]{GJ}):
\begin{equation}\label{eq:s-fraction}
  F(t, z) = \frac{1}{1 - \displaystyle\frac{z}
                    {1 - \displaystyle\frac{t z}
                    {1 - \displaystyle\frac{2t z}
                    {1 - \displaystyle\frac{2t^2 z}
                    {1 - \displaystyle\frac{3t^2 z}
                    {1 - \displaystyle\frac{3t^3 z}
                    {1 - \displaystyle\frac{4t^3 z}
                    {1 - \displaystyle\frac{4t^4 z}
                    {1 - \cdots}}}}}}}}}.
\end{equation}
Here, for $k \geq 0$, the $(2k)$th term is $(k+1)t^k z$, and the $(2k+1)$st term is $(k+1)t^{k+1} z$.

\begin{table}
\begin{center}
\begin{tabular}{cc}
\toprule
$k$ & $\abs{\set{w \in S_n : \dep(w) = k}}$ \\
\midrule
0 & 1 \\[2ex]
1 & $n-1$ \\[2ex]
2 & $\binom{n-2}{2} + 3(n-2)$ \\[2ex]
3 & $\binom{n-3}{3} + 6\binom{n-3}{2} + 9(n-3)$ \\[2ex]
4 & $\binom{n-4}{4} + 9\binom{n-4}{3} + 27\binom{n-4}{2} + 31(n-4) + 4$ \\[2ex]
5 & $\binom{n-5}{5} + 12\binom{n-5}{4} + 54\binom{n-5}{3} + 116\binom{n-5}{2} + 113(n-5) + 24$ \\[2ex]
6 & $\binom{n-6}{6} + 15\binom{n-6}{5} + 90\binom{n-6}{4} + 282\binom{n-6}{3} + 489\binom{n-6}{2} + 443(n-6) + 148$ \\[2ex]
7 & $\binom{n-7}{7} + 18\binom{n-7}{6} + 135\binom{n-7}{5} + 556\binom{n-7}{4} + 1375\binom{n-7}{3} + 2074\binom{n-7}{2} + 1809(n-7) + 744$ \\
\bottomrule
\end{tabular}
\end{center}
\caption{Formulas for the number of permutations with depth $k$.}\label{tab:depk}
\end{table}

\begin{remark}
Timothy Walsh (conveyed to us through personal communication) has, by extensive computation, given formulas for small values of depth, shown in \autoref{tab:depk}. The formulas are polynomials of degree $k$ that hold for $n\geq k$.

Each of these polynomials gives the coefficients of $t^k z^n$ in $F(t, z)$ for a fixed value of $k$, so they can be extracted from \autoref{eq:s-fraction} by computing modulo $t^{k+1}$, in which case the continued fraction eventually terminates. The result is a rational function of $t$ and $z$, and setting $t=0$ in a suitable derivative of this bivariate rational function gives a univariate rational function of $z$, say $F_k(z)$, whose power series expansion gives the desired coefficients.

It is somewhat surprising that the coefficient of $z^n$ in the power series expansion of $F_k(z)$ is polynomial in $n$; one would usually expect the coefficients to grow exponentially. This implies that the denominator of $F_k(z)$ is $(1 - z)^k$, which is not immediately obvious to us.

However, this polynomiality can be seen directly from the combinatorics of Motzkin paths with fixed area. Indeed, for any fixed area $k$, there are finitely many Motzkin paths with area $k$ that do not contain an $H$ step at height 0, and we can get an arbitrary Motzkin path by inserting $H$ steps at height 0. The number of ways of doing this is a binomial coefficient, which is polynomial in $n$.
\end{remark}

%---------------------------------------------------------------

\end{document}